\documentclass[12pt]{article}
\usepackage{amsmath,amsfonts,amssymb,theorem}

 \numberwithin{equation}{subsection}
 \numberwithin{footnote}{section}

 \newtheorem{rtheorem}[subsection]{Reduction Theorem}

 \newtheorem{thm}[subsection]{Theorem}

 \newtheorem{fga}[subsection]{The FGA Conjecture} 
 \newtheorem{ccs}[subsection]{The CCS Conjecture}

 {
 \theorembodyfont{\rmfamily}
 \newtheorem{defn}[subsection]{Definition}
 
 \newtheorem{exa-cr}[subsection]{Example--Construction}

 \newtheorem{rem}[subsection]{Remark}

 }
 \newcommand{\qed}{\ifhmode\unskip\nobreak\fi\quad\ensuremath\square}
 \newenvironment{proof}{\paragraph{Proof}}{\par\medskip}

\title{On Shokurov's Log Flips: The 3 dimensional case}

\date{Dec 2002} 
\author{Caucher Birkar\thanks{\small School of Mathematical Sciences, Nottingham University, Nottingham, UK.}}

\begin{document}
 \maketitle

\begin{abstract} I try to explain the new ideas used in the recent paper by V.V. Shokurov [Sh4] on the existence of log flips (pl flips), in the 3 dimensional case. 

 \end{abstract}

\tableofcontents
\clearpage


\clearpage{\pagestyle{empty}\clearpage}

\section{\normalsize{Introduction}}

 The birational classification of algebraic varieties in dimensions more than 2 has fundamental differences from the classification of curves and surfaces. In the case of curves, essentially we don't have birational classification because any two birational normal projective curves are isomorphic. The case of surfaces is more complicated but still we don't face much difficulties. The exceptional locus is always a bunch of rational curves and we always deal with nonsingular surfaces as far as we are concerned about the classification of nonsingular surfaces. But in the case of 3-folds or higher dimensions the exceptional locus can be of codimension more than 1 and this creates fundamental difficulties. It creates rough sorts of singularities (even non-$\mathbb{Q}$-Gorenstein), then we have to do an operation to get rid of this singularities, which is called $flip$ (see below for the definition). The flip operation proved by M. Reid [R3] for toric varieties, turned out to be extremely difficult in the general case. The general case was proved by S. Mori [M] in dimension 3 with terminal singularities. In more general settings and using quite different methods , V. Shokurov [Sh2] proved the flip problem in dimension 3 with log terminal singularities. Recently Shokurov in his fundamental paper [Sh4] created new powerful methods which are able to prove the flip problem shortly in dimension 3 and more complicated in dimension 4. The paper is very technical and still not digested by algebraic geometers. \\
This note is intended for two group of people. Those who may read this before starting Shokurov's Marathon. And those who do not want to get into technical details. \\         
 Finally I am grateful to Prof V.V. Shokurov for reading the manuscript and his valuable comments.

\section{\normalsize{Different Sorts of Flips}}

 For basic definitions I refer you to [KM]. Let (X,B) be a Klt pair and $f: X\rightarrow Z$ a birational contraction. In this section we assume that $\rho(X/Z)=1$. 

\begin{defn} $f$ is called a flipping contraction if the followings hold:

\begin{description}
 \item[Klt] $K_{X}+B$ is Klt.
 \item[small] $f$ is a small contraction i.e. $codim ~Exc(f)>1$.   
 \item[Fano] $K_{X}+B$ is $f$-antiample.
\end{description}
\end{defn}
  Now $(Z, f_{*}B)$ does not have Klt singularities (actually even it is not $\mathbb{Q}$-Gorenstein), so we have to replace it with some other varieties with Klt singularities, in hope of getting a better model for (X,B). The nominated variety is a pair $(X^{+},B^{+})$ and a map $f^{+}: X^{+}\rightarrow Z$ such that:

\begin{description}  
 \item[small] $f^{+}$ is a small contraction.   
 \item[$\mathbb{Q}$-Gorenstein] $K_{X^{+}}+B^{+}$ is $\mathbb{Q}$-Cartier.
 \item[Compatible] $B^{+}$ is the strict transform of $B$.
 \item[Ample] $K_{X^{+}}+B^{+}$ is $f^{+}$-ample. 
\end{description}

 I haven't assumed much about the singularities of $(X^{+},B^{+})$, but it turns out that its singularities are at least as good as $(X,B)$. This process is called $K_{X}+B$-flip if it exists. Now one might ask a stupid question: is this really the only choice we have to do? At least it looks very natural, because we make the log canonical divisor more nef and this is the ultimate goal of the whole program. \\

 It is well known that this problem is equivalent to the finite generation of the following  sheaf of graded $\mathcal{O}_{Z}$-algebras:

 \[\mathcal{R}= \mathcal{R}_{X/Z}(K_{X}+B)=\mathcal{R}(X/Z,K_{X}+B)=\bigoplus_{i=0}^{\infty}f_{*}\mathcal{O}_{X}(i(K_{X}+B))\] 

  If this algebra is fintely generated then we take $X^{+}=Proj \mathcal{R}$. This is the first step toward the algebraisation of the problem. Algebraic methods are usually much more powerful and better in higher dimensions. Shokurov's idea here is to reduce the problem to lower dimensions, that is to use induction. So he reduces the problem to a special kind of flips, pl flips where the reduced part of the boundary is not zero. This enables him to use adjunction and good properties of components in the reduced part and then to restrict the above algebra to the intersection of these components. He proves the existence of more general kind of flips than what I defined above, that is log flips. Now I give the definition of log and pl flips:       

\begin{defn}
  $f: X\rightarrow Z$ a birational contraction is called a log flipping contraction if:
\begin{description}
 \item[Klt] $K_{X}+B$ is Klt.
 \item[Antinef] $K_{X}+B$ is $f$-antinef.
 \item[Small] f is small.
\end{description}
\end{defn}

 And the log flip operation or $K_{X}+B$-flip is similar to what I defined above, replacing $f^{+}$-ample by $f^{+}$-nef for $K_{X^{+}}+B^{+}$.\\ 

\begin{defn}
Let $1\leq s$  and $S=\sum_{i=1}^{s}S_{i}$ be a sum of reduced Weil divisors on $X$ then $f: X\rightarrow Z$ a birational contraction is called a pl contraction if:
\begin{itemize}
 \item $K_{X}+B+S$ is dlt and plt if $s=1$.   
 \item $K_{X}+B+S$ is $f$-antiample.
 \item each $S_{i}$ is $\mathbb{Q}$-Cartier and $S_{i}\sim_{\mathbb{Q}}r_{i,j}S_{j}$ for rational numbers $r_{i,j}>0$.
\end{itemize}
 We say $f$ is an elementary pl contraction if in addition the followings hold:
   
\begin{itemize}
 \item $f$ is extremal, that is the relative Picard number $\rho(X/Z)=1$.   
 \item $S$ is $f$-antiample.
 \item $f$ is small.
 \item $X$ is $\mathbb{Q}$-factorial and projective/Z.
\end{itemize}
\end{defn}

 The $S$-flip for this contraction is called pl flip if it exists. And  $S$-flip is as in definition 2.1 replacing $K_{X}+B$ and $K_{X^{+}}+B^{+}$ by $S$ and $S^{+}$ respectively.
 
\begin{rem}
Note that I have defined three different sort of flips: first one is flip which is the traditional one, the second is log flip which is the generalised sort of flip and the third one is pl flip which is an auxiliary tool to construct log flips.
\end{rem}

The following theorem shows why we are interested in pl flips.

\begin{rtheorem}
Log flips exist in dimension $n$ if the followings hold:
\begin{description}
 \item[Pl Flips] Pl flips exist in dimension $n$.
 \item[Special Termination] special termination holds in dimension $n$. 
\end{description}
\end{rtheorem}
\begin{proof}
See Shokurov [Sh4]. The main idea is to choose a good reduced Cartier divisor $H$ on $Z$ such that it contains all singularities of $Z$ and singularities of the push down of the boundary on $Z$. And moreover that the components of $^{*}H$ on any model $W$ of $Z$ generate the Neron-Severi group of $W$. Then we take $R\rightarrow X $ a log resolution for the pair $(X,B)$ and put $D=B^{-}+H^{-}+\sum E_{i}$ where the superscript $^{-}$ stands for the strict birational trasform and $E_{i}$ are all exceptional divisors of the resolution. Now we start running the LMMP for the pair $(R,D)$. In each step, discarding the relatively ample components of $D$, we face a pl contraction or a divisorial contraction if we choose our contractions to be extremal. So we are fine by the assumptions on pl flips. The program terminates by assumption on special termination.  
\end{proof}
\begin{rem}
Special termination claims that if we have a sequence of flips then after a finite number of steps the flipping locus does not intersect the reduced part of the boundary. More generally it does not intersect any log canonical centre on $X$. In the proof above in each step some component in the reduced part of the boundary is relatively negative (by assumptions on $H$) so the special termination applies to this case.
\end{rem}  
\begin{rem}
 In order LMMP in lower dimensions implies the special termination in dimension $n$ [Sh4, 2.3], so we don't have to worry about the special termination in the 4 dimensional case. The important thing is to prove the existence of pl flips.\\ To prove the speial termination using LMMP in lower dimensions, we note that log canonical centres of $(R,D)$ (in particular irreducible components of the reduced part of $D$) on $R$ are located in the local intersection of irreducible components of the reduced part of $D$. If $(R_{i},B_{i})\rightarrow (R_{i+1},B_{i+1})/T_{i}$ is a sequence of flips and $\omega$ a log canonical centre where $R_{0}=R$. Then  using adjunction we get a sequence of birational operations $(\omega_{i},B_{\omega_{i}})\rightarrow \omega_{i+1},B_{\omega_{i+1}})/\gamma_{i}$ where we may have both divisorial and small contractions /$\gamma_{i}$. We have to get rid of the divisorial ones (using versions of $difficulty$ introduced by Shokurov) and get a sequence of log flips for 
$(\omega,B_{\omega})$ and use the LMMP to conclude that the original sequence of flips induces isomorphisms on $(\omega,B_{\omega})$. Then we can easily get the special termination from this. Because if any flipping curve $C_{i}$ intersects $\omega_{i}$ then $C_{i}.\omega_{i}>0$ (note by the above $\omega$ does not contain any flipping curve ). So we have $C_{i+1}. \omega_{i+1}<0$ for some flipped curve $C_{i+1}$, that is, $\omega_{i+1}$ contains $C_{i+1}$ which is a contradiction. 

\end{rem}

\section{\normalsize{Reduction to Lower Dimensions and b-divisors}}

To prove the existence of pl flips now we know how to use induction (of course after Shokurov!). The targeted lower dimensional variety is the intersection of all $S_{i}$ given in the definition of pl flips. $Y=\bigcap_{i=1}^{s}S_{i}$ is called the core of $f$ and its dimension, $d$, is called the core dimension. The smaller is $d$ the easier is life. $Y$ is normal by the fact that $(X, K_{X}+B+S)$ is a dlt pair. It is irreducible near the fibers of a point $P$ on $Z$ (we can shrink Z as our problem is local with respect to $Z$). Using adjunction we also know that the new pair $(Y/T,B_{Y})$ is Klt where $T=f(Y)$. \

Moreover the special termination is proved up to dimension $4$, so the existence of pl flips in dimension $4$ implies the existence of all log flips in dimension $4$.\\
 The existence of pl flips is also equivalent to the finite generation of a graded sheaf of algebras, that is $\mathcal{R}_{X/Z}(D)$ for a suitable $D\sim_{\mathbb{Q}}S$. Now we can restrict this algebra naturally to $Y$ via maps $r_{i}: \mathcal{O}_{X}(iD)\rightarrow \mathcal{O}_{Y}(iD\vert_{Y})$ and denote it  by $\mathcal{R}\vert_{Y}$. Unfortunately the resulting algebra is not divisorial i.e. it is not as $\mathcal{R}\vert_{Y}=\bigoplus_{i=0}^{\infty}f_{*}\mathcal{O}_{Y}(i(D))$ for a divisor $D$ on $Y$. Now the beautiful idea of Shokurov remedies this difficulty, the notion of b-divisor or birational divisor. An extention ${\mathcal{\overline{R}}\vert_{Y}}$, of the above algebra, is pseudo b-divisorial algebra. i.e. there are b-divisors $\mathcal{M}_{i}$ such that 

\[\mathcal{\overline{R}\vert}_{Y}=\bigoplus_{i=0}^{\infty}f_{*}\mathcal{O}_{Y}(\mathcal{M}_{i})\]   
The finite generation of these algebras are equivalent [Sh4, 3.43 and 4.15].
A b-divisor $\mathcal{D}$ over $Y$ is a formal sum of prime divisors i.e. $\mathcal{D}=\sum_{i=1}^{\infty}d_{i}D_{i}$ , with integral coefficients, where $D_{i}$ are on birational models of $Y$ such that there are just finitely many of them on any single birational model. And also it should be compatible with push down of divisors. $\mathcal{D}_{W}=\sum_{i=1}^{\infty}d_{i}D_{i}$ , such that $D_{i}$ is a divisor on $W$, is the trace of $\mathcal{D}$ on $W$ where $W$ is a birational model of $Y$. Let $D$ be a Cartier divisor on $Y$ then $\overline{D}$ is a b-divisor which has the trace $^{*}D$ on any model of $Y$ over $Y$. This b-divisor is called the Cartier closure of $D$. Similarly b-divisors are defined over the fields of rational and real numbers. Actually $\mathcal{M}_{i}$ above are defined as \[\mathcal{M}_{i}= lim sup\{-\overline{(s)}: s\in \mathcal{R}_{i}\}\] where $\mathcal{\overline{R}\vert}_{Y}=\bigoplus_{i=0}^{\infty}\mathcal{R}_{i}$. 

\begin{thm}
$\mathcal{R}$ is f.g. if and only if $\mathcal{R} \vert_{Y}$ is f.g. if and only if $\mathcal{\overline{R}\vert}_{Y}$ is f.g. 
\end{thm}
\begin{proof}
See [Sh4, 3.43] and [Sh4, 4.15]. 
\end{proof}

And now we convert our problem to another about b-divisors. Put $\mathcal{D}_{i}=\mathcal{M}_{i}/i$.

\begin{thm}[Limiting Criterion]
$\mathcal{\overline{R}}\vert_{Y}$ is f.g. if and only if the system $\{\mathcal{D}_{i}\}_{i=0}^{\infty}$ stabilises i.e. $\mathcal{D}_{i}=\mathcal{D}$ for all large $i$ where $\mathcal{D}=lim_{i \rightarrow \infty} \mathcal{D}_{i}$. 
\end{thm}

\begin{proof}
See [Sh4, 4.28]. 
\end{proof}

In practise, first we always try to prove that $\mathcal{D}$ is a b-divisor over $\mathbb{Q}$ and then prove that the system stabilises. To prove this rationality condition Shokurov introduced the notion of (asymptotic) saturation of linear systems and proved that our system has that property [Sh4, section 4]. It means:

And now one of the central notions in [Sh4]:
\begin{description}
 \item[log canonical asymptotic saturation] There is an integer $I$, called the index of the saturation, such that for any i and j which satisfies I$\vert$ i,j we have 

\[Mov ~\ulcorner j \mathcal{D}_{i}+\mathcal{A} \urcorner_{W} \leq (j \mathcal{D}_{j})_{W}\]

 for any high resolution $W$, where $\mathcal{A}$ is the discrepancy b-divisor i.e. $\mathcal{A}_{W}=K_{W}-{^{*}(K_{Y}+B_{Y})}$.  
\end{description}

\begin{rem}[Truncation Principle] It is easy to prove that $\bigoplus_{i=0}^{\infty}\mathcal{R}_{i}$ is f.g. if and only if $\bigoplus_{i=0}^{\infty}\mathcal{R}_{il}$ is f.g. for a natural number $l$. So we always may replace our sequence $\mathcal{D}_{i}$ with  $\mathcal{D}_{il}$ without mentioning it.
\end{rem}

\section{\normalsize{The FGA Conjecture}}

Shokurov considers much more general settings and proposes the following conjecture which implies our original conjecture and he proves this in the two dimensional case:

\begin{fga}
Let $(Y/T, B)$  be a weak log Fano contraction (in particular Klt) then any system of b-divisors $\{\mathcal{D}_{i}\}_{1}^{\infty}$ which satisfies the followings, stabilises.
\begin{itemize}
 \item $_{*}\mathcal{O}_{Y}(i\mathcal{D}_{i})$ is a coherent sheaf on $T$ for all $i$. 
 \item log canonical asymptotic sturation. 
 \item convexity i.e. $i\mathcal{D}_{i}+j\mathcal{D}_{j}\leq (i+j)\mathcal{D}_{i+j}$ for all $i,j$.
\end{itemize} 
\end{fga}

Now let have a look at the one dimensional case of this conjecture. For curves, b-divisors are usual divisors. Let $(C/pt. ,B)$ be a klt pair, then the saturation looks like the following, componentwise:

\[ \ulcorner j{d}_{i}+{a} \urcorner \leq j{d}_{j}\]

where $a=-b$, $b<1$ and note that divisors with high degree have no fixed part, so are movable. By definition $d=\lim_{i\rightarrow \infty}d_{i}$ so we have:
\[ \ulcorner j{d}-{b} \urcorner \leq j{d}\]
This last formula implies that $d$ is a rational number. If this is not the case the set $\{<jd>: j\in \mathbb{N}\}$ is dense in the interval $[0,1]$, where $<x>$ stands for the fractional part of $x$. This fact and the fact that $b<1$ implies that $d$ should be rational. So for some $j$ we have $ j{d}+\ulcorner -b \urcorner \leq jd_{j}$ so $d \leq d_{j}$ and then $d=d_{j}$ (for infinitely many $j$).
This approximation procedure is an essential part of the problem also in higher dimensions. To reach this, we use the fact that semiampleness is an open condition. But we should first prove that $\mathcal{D}_{Y}$ is semiample on certain models. We can make $\mathcal{D}_{Y}$ nef using LMMP and then the weak log Fano condition plays its role: in this case, nef divisors are semiample (/$T$). \\

\section{\normalsize{Finding Good Models}}

B-divisors $i\mathcal{D}_{i}$ appeared in section 3 in the restriction algebra have many good properties. I mentioned the log canonical asymtotic saturation property. $i\mathcal{D}_{i}$ are b-free (and so b-nef) which means that there is a model $W$ such that $i\mathcal{D}_{i}=\overline{i\mathcal{D}_{i}}_{W}$ and  ${i\mathcal{D}_{i}}_{W}$ is a free divisor and in particular a nef divisor, but for different $\mathcal{D}_{i}$ we have different $W$ (the ultimate goal is to prove that many of them share the same model). $\mathcal{D}_{i}$ also have all properties sorted in the statement of the FGA conjecture. In the last section I mentioned the approximation procedure and used it to prove the FGA conjecture in one dimensional case. But in higher dimensions it is more complicated. To use this method we can first make infinitely many $\mathcal{D}_{i}$ nef over a single model $/T$ of $(Y,B)$ and still not to loose the weak log Fano condition and other mentioned properties of $\mathcal{D}_{i}$. We replace $(Y,B)$ by this model and again show by $(Y,B)$. In the course of obtaining this model the boundary may increase (see [Sh4, Example 5.27] for full details). Now all ${\mathcal{D}_{i}}_{Y}$ being nef$/T$  implies that $\mathcal{D}_{Y}$ is also nef$/T$ and so semiample$/T$. The dificulty is that we do not know if  $D=\mathcal{D}_{Y}$ is a $\mathbb{Q}$-divisor so we can not simply say that a multiple of it, is free$/T$. But we know that $\mathbb{Q}$-divisors very close to it are eventually free.Assuming that $D$ is not a $\mathbb{Q}$-divisor and  using Diophantine approximation we can get  $\mathbb{Q}$-divisors  $\{D_{\alpha}\}_{\alpha \in \mathbb{N}}$ such that:
\begin{itemize}
 \item $D_{\alpha}\nleq D$ for any $\alpha$. 
 \item $\alpha D_{\alpha}$ is free.
 \item for any $\epsilon$ there is an $N$ such that $|\alpha D_{\alpha}-\alpha D|\leq \epsilon$ if $N\leq\alpha$.
\end{itemize}
To prove our stabilisation it is enough to prove that for a crepant model $(U,B_{U})/T$ of $(Y,B)/T$ we have the followings:

\begin{enumerate}
 \item $\mathcal{D}_{i}=\overline{(\mathcal{D}_{i})}_{U}$.
 \item $\mathcal{D}_{U}=(\mathcal{D}_{i})_{U}$. 
for infinitely many $i$.
\end{enumerate}
 
 Now lets consider how we can solve this problem using informations above assuming that the approximation has been caired out on $(Y',B')$ a crepant model of $(Y,B)/T$:
  
\[Mov ~\ulcorner j \mathcal{D}_{i}+\mathcal{A} \urcorner_{W}=\] \[
Mov ~\ulcorner (\mathcal{A}-j{\overline{D}_{i}}_{Y'}+j\mathcal{D}_{i})+(j{\overline{D}_{i}}_{Y}-j\overline{D})+(j\overline{D}-j\overline{D_{\alpha}})+j\overline{D_{\alpha}} \urcorner_{W} \leq (j \mathcal{D}_{j})_{W} \leq (j\mathcal{D})_{W}\]\\

 Denote $\overline{\mathcal{N}}_{Y'} - \mathcal{N}$ by $\mathcal{E}_{\mathcal{N},Y'}$ for any b-divisor $\mathcal{N}$.  By negativity lemma $0 \leq \overline{\mathcal{N}_{Y'}}-\mathcal{N}$ if $\mathcal{N}$ is b-nef.  To get a contradiction, the only bothering term in the above formula is $-j({\overline{\mathcal{D}}_{i}}_{Y'}-\mathcal{D}_{i})+\mathcal{A}=-j\mathcal{E}_{\mathcal{D}_{i},Y'}+\mathcal{A}$. If we can prove that $0 \leq -r_{i}\mathcal{E}_{\mathcal{D}_{i},Y'}+\mathcal{A}$ for all $i$ such that $r_{i}\rightarrow \infty$, then we get a contradiction, even if we prove that $ 0\leq \ulcorner -r_{i}\mathcal{E}_{\mathcal{D}_{i}}+\mathcal{A} \urcorner$. This inequality is one of the most important things that Shokurov tries to prove and this leads to the CCS conjecture. Shortly, it is important to have the followings on a crepant model $(Y',B')/T$ of $(Y,B)/T$ which is called a prediction model for the problem:

\begin{description}
 \item[semiampleness] $\mathcal{D}_{Y'}$ is semiample.   
 \item[canonical asymptotic confinement] There are positive real numbers $r_{i}$ such that $r_{i} \rightarrow \infty$ and $0\leq -r_{i}\mathcal{E}_{\mathcal{D}_{i},Y'}+\mathcal{A}$ holds for any $i$.
\end{description}
I discussed the first condition and the second one will be discussed in the next section. Note that the second one implies that $0 \leq\mathcal{A}(Y',B')$ since $\mathcal{D}_{i}$ are b-nef. This means that $(Y',B')$ is canonical. So this predicts that to get the asymptotic confinement it is better to work on a terminal crepant model of $(Y,B)/T$ as we know there exists.
\\

\section{\normalsize{The CCS Conjecture}}

The conditions at the end of the last section are sufficient to solve our problem but the second is not necessary at all. On the other hand this condition has one very important advantage, that is $\mathcal{E}_{N}=\mathcal{E}_{\acute{\mathcal{N}}}$ if $\mathcal{N}\sim {\mathcal{\acute{N}}}$ or even if $\mathcal{N}\equiv {\mathcal{\acute{N}}}$. And we know that the saturation condition is preserved under linear equivalent changes. So we may move our divisors linearly and use their freeness properties.  
 
\begin{description}
 \item[canonical confinement of singularities(CCS)] Let $\{{D}_{\alpha}\}_{\alpha \in A}$ be a set of divisors on $(Y',B_{Y'})$. We say that singularities of this divisors is confined up to linear equivalence if there is $0<c$ such that for any $\alpha$ there is ${D'}_{\alpha}\in |D_{\alpha}|$ s.t. the pair $(Y',B_{Y'}+c{D'}_{\alpha})$ is canonical.  
\end{description}

 Note that the general member of a free linear system is reduced and irreducible. So any such divisor can be confined by a $c$ which just depends on the model and not on the free divisor. Actually all b-divisors in the conjecture bellow  are free on the terminal model but a bounded family of them. The bounded family in the conjecture corresponds to that bounded family of divisors, because each of those divisors is free on a model depending on the divisor. \\
Back to our set of b-divisors $\{\mathcal{D}_{i}\}$. Suppose there is a $c$ which confines the singularities of these divisors on $(Y',B_{Y'})$ and ${{\mathcal{D'}}_{i}}_{Y'}\in |{\mathcal{D}_{i}}_{Y'}|$ (see 6.1). Then for any model $W$ over $Y'$ we have:

\[(\mathcal{A}(Y',B_{Y'})-ic\mathcal{E}_{\mathcal{D}_{i},Y'})_{W}=(\mathcal{A}(Y',B_{Y'})-ic\mathcal{E}_{{\mathcal{D'}_{i}},Y'})_{W}=\] \[K_{W}-^{*}(K_{Y'}+B_{Y'})-ic(\mathcal{E}_{{\mathcal{D'}_{i}},Y'})_{W}=K_{W}-^{*}(K_{Y'}+B_{Y'}+ci{{\mathcal{D'}_{i}}_{Y'}})+ci({\mathcal{D'}_{i}})_{W}\]\[\geq ci({\mathcal{D'}_{i}})_{W}\geq 0.\]

In other words, we have the asymptotic confinement for the b-divisors $\mathcal{D}_{i}$ over $(Y',B_{Y'})$. Now, can we find such a model? This is what CCS conjecture is about.


Roughly speaking the CCS conjecture is as follows:   

\begin{ccs}
Let $\mathfrak{M}(Y',B_{Y'})$ denote the set of b-free b-divisors which are  log canonically saturated (i.e. $Mov \ulcorner \mathcal{M}+\mathcal{A} \urcorner \leq \mathcal{M}$. Then there is a bounded family of models  on which $\mathfrak{M}(Y',B_{Y'})$ has canonically confined singularities. In other words, there is $c>0$ s.t. for each $\mathcal{M}\in\mathfrak{M}(Y',B_{Y'})$ there is a creap terminal model $(Y'_{\mathcal{M}},B_{\mathcal{M}})$ and ${\mathcal{M'}}\in |\mathcal{M}|$ such that $B_{Y_{\mathcal{M}}}+c{\mathcal{M'}_{Y'_{\mathcal{M}}}}$ is canonical. Moreover if $(Y,B_{Y})/T$ is birational then this family can be taken finite (And for our problem it can be taken just one model). 
\end{ccs}   

See [Sh4, 6.14]. Now by asymptotic saturation for $\{\mathcal{D}_{i}\}$ we get the canonical saturation for $\mathcal{M}_{j}=j\mathcal{D}_{j}$ :

\[Mov ~\ulcorner j \mathcal{D}_{j}+\mathcal{A} \urcorner_{W} \leq (j \mathcal{D}_{j})_{W}\]

where we take $i=j$. So we may apply the above conjecture to prove the canonical confinement of singularities. Also on this model we have the semiampleness property because it is a crepant model of $(Y,B)$, so this gives a solution to our problem. This conjecture has been proved up to dimension 2 [Sh4, 6.25 and 6.26].

\section{\normalsize{References:}}
\begin{description}

 \item[[A1]] ~F. Ambro; ${Shokurov's ~pl ~flips}$, preprint, cam.dpmms.cam.ac.uk/~fa239/pl.ps
 \item[[A2]] ~F. Ambro; ${On ~minimal ~log ~discrepancies}$; Math. Res. Letters 6(1999), 573-580. 
 \item[[F]] ~ O. Fujino; ${Private ~ notes ~ on ~ special ~termination ~ and ~ the ~reduction~ theorem}$. 
 \item[[H]] ~R. Hartshorne; ${Algebraic ~Geometry}$, Springer-Verlag, 1977. 
 \item[[KMM]] ~Y. Kawamata, K. Matsuda, K. Matsuki; ${Introduction ~to ~the ~minimal}$\\${ ~model ~problem}$, in Algebraic Geometry (Sendai, 1985) Adv. Stu. Pure Math. 10 (1987), kinokuniya, 283-380.
 \item[[K1]] ~J. Kollar; ${Singularities ~of ~pairs}$, preprint, 1996. 
 \item[[KM]] ~J. Kollar, S. Mori; ${Birational ~geometry ~of ~algebraic ~varieties}$, Cambridge University Press, 1998.
 \item[[M]] S. Mori, ${Flip ~theorem ~and ~the ~existence ~of ~minimal ~model ~for ~3-folds}$, J. AMS 1 (1988) 117-253.  
 \item[[R1]] ~M. Reid; ${Young ~person's ~guide ~to ~canonical ~singularities}$, Algebraic geometry, Bowdoin, 1985 (Brunswick, Maine, 1985), 345--414, Proc. Sympos. Pure Math., 46, Part 1, Amer. Math. Soc., Providence, RI, 1987. 
 
 \item[[R2]] ~M. Reid; ${Chapters ~on ~algebraic ~surfaces}$, in 1AS/Park City Mathematics Series 3 (1997) 5-159.
 \item[[R3]] ~M. Reid; ${Decomposition ~of ~toric ~morphisms}$, in Arithmetic and Geometry, papers dedicated to I.R. Shafarevich, Birkhauser 1983, Vol II, 395-418.  

 \item[[Sh1]] ~V.V. Shokurov; ${The ~nonvanishing ~theorem}$, Math. USSR  Izvestija,26(1986) 591-604. 
 \item[[Sh2]] ~V.V. Shokurov; ${3-fold ~log ~flips}$, Russian Acad. Sci. Izv. Math., 40 (1993) 95-202. 

 \item[[Sh3]] ~V.V. Shokurov; ${3-fold ~log ~models}$, Algebraic geometry, 4. J. Math. Sci. 81 (1996), no. 3, 2667--2699.. 
 \item[[Sh4]] ~V.V. Shokurov; ${Pl ~flips}$, Proc. Steklov Inst. v. 240, 2003.
 \item[[Sh5]] ~V.V. Shokurov; ${Letters ~of ~a ~birationalist ~IV: ~Geometry ~of ~log ~flips}$, Alg. Geom. A volume in memory
of Paolo Francia. Gruyter 2002, 313-328.
 \item[[T]] ~H. Takagi; ${3-fold ~log ~flips ~according ~to ~V. ~Shokurov}$, preprint 1999.

\end{description}
   
\end{document}